\documentclass[amscd,amssymb,11pt,reqno]{amsart}

\newtheorem{Thm}{Theorem}[section]

\newtheorem{Lem}[Thm]{Lemma}

\usepackage{graphicx}
\usepackage{latexsym}
\begin{document}


\vspace{1.5 cm}

\title[Non-intersection bodies.]
      {Non-intersection bodies all of whose central sections are intersection bodies.}\

\author{M.Yaskina}

\address{M.Yaskina, Department of Mathematics, University of Missouri, Columbia, MO 65211, USA}
\email{yaskinam@math.missouri.edu}

\begin{abstract}
We construct symmetric convex bodies that are not intersection bodies, but all of their central
hyperplane sections are intersection bodies. This result extends the studies by Weil in the case of
zonoids and by Neyman in the case of subspaces of $L_p$.
\end{abstract}

\subjclass[2000]{52A20, 52A21, 46B20.}

\maketitle

\section{Introduction}
The concept of an intersection body was introduced by Lutwak \cite{Lu} in 1988. Let $K$ and $L$ be
origin symmetric star bodies in $\mathbb{R}^n$. Following \cite{Lu}, we say that $K$ is the {\it
intersection body of} $L$ if the radius of $K$ in every direction is equal to the volume of the
central hyperplane section of $L$ perpendicular to this direction, i.e. for every $\xi \in
S^{n-1}$,
$$\|\xi\|_K^{-1}=\mathrm{vol}_{n-1}(L\cap\xi^\bot).$$
The closure in the radial metric of the class of intersection bodies of star bodies gives the class
of {\it intersection bodies}.

Intersection bodies played an important role in the solution of the Buse\-mann-Petty problem (see
\cite{GKS} and \cite{Zh2} for the solution and historical details). Posed in 1956, \cite{BP}, the
Busemann-Petty problem asks the following. Let $K$ and $L$ be two origin-symmetric convex bodies in
$\mathbb{R}^n$ so that the $(n-1)$-dimensional volume of every central hyperplane section of $K$ is
smaller that the same for $L$. Does it follow that the $n$-dimensional volume of $K$ is smaller
than the $n$-dimensional volume of $L$? The answer turns out to be affirmative for dimensions $n\le
4$ and negative for $n\ge 5$.

The connection between intersection bodies and the Busemann-Petty problem was found by Lutwak
\cite{Lu}. First, the answer to the problem is affirmative if $K$ is an intersection body and $L$
is any origin-symmetric star body. On the other hand, if $L$ is an origin-symmetric convex body
that is not an intersection body, one can perturb $L$ to construct a counterexample to the
Busemann-Petty problem. Hence, a solution to the problem in $\mathbb{R}^n$ is affirmative if and
only if every infinitely smooth origin-symmetric convex body in $\mathbb{R}^n$ is an intersection
body, which is the case for dimensions $n\le 4$.
Examples of non-intersection bodies in dimensions $5$ and higher were constructed in \cite{Ga1},
\cite{Zh1}, \cite{GKS}, \cite{K1.5}.

In this paper we are interested in the following problem. Does there exist a convex
body $K$ that is not an intersection body, but every its section by a central
hyperplane is an intersection body? We construct an example of such a body for
dimensions $n\ge 5$. Our result can also be considered as a new way of constructing
non-intersection bodies.

This paper was motivated by results of W.Weil \cite{W} and A.Neyman \cite{N}. In 1982 W.Weil,
\cite{W}, showed that it is not possible to characterize zonoids by means of their projections. He
constructed a convex body in $\mathbb{R}^n$ ($n \ge 3$) that is not a zonoid but all its
projections onto hyperplanes are zonoids. A.Neyman in \cite{N} showed that there are
$n$-dimensional normed spaces that do not embed in $L_p$, but all their $(n-1)$-dimensional
subspaces embed in $L_p$ for $p>0$. He used this to prove that for $p>0$, $p\neq 2$, $L_p$ is not
characterized by a finite number of equations. Let us note that M.Burger in \cite{B} used another
approach to show that, for $n \ge 3$, zonoids cannot be characterized by a finite number of
piecewise inequalities. A.Koldobsky in \cite{K2} introduced the concept of embedding of a normed
spaces in $L_p$, $p<0$, and proved that intersection bodies are the unit balls of spaces that embed
in $L_{-1}$. Therefore, our result can be considered as an extension of Neyman's example to
negative $p$.



\section{Main Results}
Our main tool is the Fourier transform of distributions. The Fourier transform of a distribution
$f$ is defined by $\langle\hat{f}, \phi\rangle= \langle f, \hat{\phi} \rangle$ for every test
function $\phi$ from the Schwartz space $ \mathcal{S}$ of rapidly decreasing infinitely
differentiable functions on $\mathbb R^n$. For any even distribution $f$, we have $(\hat{f})^\wedge
= (2\pi)^n f$.

A distribution is {\it positive definite} if its Fourier transform is a positive distribution in
the sense that $\langle \hat{f},\phi \rangle \ge 0$ for every non-negative test function $\phi$;
see, for example, \cite[p.152]{GV}.

Let $K$ be a convex origin-symmetric body in $\mathbb R^n.$  Our definition of a convex body
assumes that the origin is an interior point of $K$ .
The {\it radial function} of $K$ is given by
$$\rho_K(x)=\max \{a>0: ax \in K \}, \ \ \ x\in \mathbb R^n \setminus \{0\} $$
The Minkowski {\it norm} of $K$ is defined by $$\|x\|_K=\min \{a\ge 0: x \in a K \}.$$ Clearly
$\rho_K(x)=\|x\|_K^{-1}$.



The main result of this paper is the following
\begin{Thm}\label{Thm:main}
There exists a convex body $K$ in $\mathbb{R}^n$, $n\ge5$, that is not an intersection body, but
for every $(n-1)$-dimensional subspace $V$ of $\mathbb{R}^n$, $K\cap V$ is an intersection body.
\end{Thm}

To construct an example of such a body we use a connection between the Fourier transform and
intersection bodies. A.Koldobsky in \cite{K1} proved that an origin-symmetric star body $K$ in
$\mathbb{R}^n$ is an intersection body if and only if $\|x\|_K^{-1}$ is a positive definite
distribution.


We will use Lemma 3.16 from \cite{K3}. It states the following:

\begin{Lem}\label{Lem:Kold}
Let $k \in \mathbb{N}\cap\{0\}$ and $f\in C^{2k}(S^{n-1})$, $f$ is even, $q\le 2k$, $q$ is not an
odd integer. Then:

\noindent (i)  The Fourirer transform of the distribution $f(\theta)r^{-n+q+1}$ is a homogeneous of
degree $-1-q$ continuous on $\mathbb{R}^n\setminus\{0\}$ function. If $q<2k$ then for every $x \in
\mathbb{R}^n$,
$$|x|_2^{2k}\left(f(\theta)r^{-n+q+1}\right)^\wedge(x)=\frac{(-1)^k\pi}{-2\Gamma(2k-q)\sin(\pi(2k-q-1)
/2)}$$ $$\times
\int_{S^{n-1}}|(x,\xi)|^{2k-q-1}\Delta^k\left(f(\theta)r^{-n+q+1})\right)(\xi)d\xi.$$ If $q=2k$
then $$|x|_2^{2k}\left(f(\theta)r^{-n+q+1}\right)^\wedge(x)=(-1)^k\pi|x|_2^{-1} $$
$$\times \int_{S^{n-1}\cap(x/|x|_2)^\bot}\Delta^k\left(f(\theta)r^{-n+q+1})\right)(\xi)d\xi, $$
where $\Delta$ is the Laplace operator in $\mathbb{R}^n$.

\noindent (ii) If $f \in C^{\infty}(S^{n-1})$ then there exist an even function $g \in
C^{\infty}(S^{n-1})$ so that for every $x=t\xi \in \mathbb{R}^n$, $t\neq 0$, $\xi \in S^{n-1}$,
$$\left(f(\theta)r^{-n+q+1}\right)^\wedge(x)= g(\xi)t^{-1-q},$$
so the Fourier transform of $f(\theta)r^{-n+q+1}$ is an infinitely smooth function on
$\mathbb{R}^n\setminus \{0\}$.
\end{Lem}

Fix a point $x_0$ on the unit sphere $S^{n-1}$. Define a function $f_\varepsilon$ as
follows:

$$ f_{\varepsilon}(x)=\left\{
\begin{array}{lll}
2 \displaystyle{e^{-\frac{|x-x_0|^2}{\varepsilon^2-|x-x_0|^2}}}& \textrm{if } |x-x_0|<\varepsilon \\
2 \displaystyle{e^{-\frac{|x+x_0|^2}{\varepsilon^2-|x+x_0|^2}}}& \textrm{if } |x+x_0|<\varepsilon \\
0 & \textrm{otherwise. }\end{array} \right. $$

Clearly $f_\varepsilon$ is an infinitely differentiable function.
Define a body $K$ by
\begin{equation}\label{Def:K}
\|x\|_K^{-1}=\left( (1- f_\varepsilon(\theta))r^{-n+1}\right)^\wedge(x), \quad x \in \mathbb{R}^n
\setminus\{0\},
\end{equation}
where $f_\varepsilon(\theta)r^{p}=f_\varepsilon\left(\displaystyle\frac{x}{|x|_2}\right)
|x|_2^{-p}$, $x=(r,\theta)$ are polar coordinates in $\mathbb{R}^n$.

The function $\|x\|_K^{-1}$ is infinitely smooth on $\mathbb{R}^n\setminus \{0\}$ by Lemma
\ref{Lem:Kold}.
 It will be shown in Lemma \ref{Lem:est} that $\left( f_\varepsilon(\theta)r^{-n+1}
\right)^\wedge(x)$ is of the order $\varepsilon^{n-2}$ uniformly with respect to $x \in S^{n-1}$,
therefore $\|x\|_K$ is positive for a small $\varepsilon$.

\begin{Lem}\label{K-int}
For any $\varepsilon>0$, $K$ is not an intersection body.
\end{Lem}

\noindent{\bf Proof.} Since $f_\varepsilon$ is an even function, we have
$$\left(\|x\|_K^{-1}\right)^\wedge=(2\pi)^n (1- f_\varepsilon(\theta))r^{-n+1},$$
which is negative for $\theta$ in some neighborhood of $x_0$. Therefore, by \cite[Theorem 1]{K1},
$K$ is not an intersection body.

\qed

\begin{Lem}\label{Lem:est}
There exist constants $D_1$, $D_2$, $D_3$ so that for every $x \in S^{n-1}$

$$\left|\left( f_\varepsilon(\theta) r^{-n+1} \right)^\wedge(x) \right| \le D_1 \varepsilon^{n-2} $$

$$\left|\displaystyle\frac{\partial}{\partial x_i}\left( f_\varepsilon(\theta)
r^{-n+1} \right)^\wedge(x) \right| \le D_2 |\ln\varepsilon|\cdot\varepsilon^{n-3}$$ and
$$\left| \displaystyle\frac{\partial^2}{\partial x_i \partial x_j}\left( f_\varepsilon(\theta)
 r^{-n+1}\right)^\wedge(x)\right| \le D_3 \varepsilon^{n-4}.$$

\end{Lem}
\noindent{\bf Proof.} First we show the estimate for the second derivative. Using the connection
between the Fourier transform and differentiation and Lemma \ref{Lem:Kold}, we get
\begin{eqnarray*}
\frac{\partial^2}{\partial x_i x_j}\left( f_\varepsilon(\theta)r^{-n+1}\right)^\wedge(x)&=&-\left(
f_\varepsilon\left(\frac{y}{|y|}\right) y_i y_j |y|^{-n+1}  \right)^\wedge (x) \\
&=& -\left( f_\varepsilon\left(\frac{y}{|y|}\right) \frac{y_i y_j}{|y|^2} |y|^{-n+3}  \right)^\wedge (x) \\
&=& -\left( g_\varepsilon(\theta) |y|^{-n+3}  \right)^\wedge (x) \\
&=&\pi|x|_2^{-3}\int_{S^{n-1}\cap (x/|x|_2)^\bot}\triangle \left(g_\varepsilon(\theta) |y|^{-n+3}
\right)(\xi) d\xi,
\end{eqnarray*}
where $$g_\varepsilon(\theta)=f_\varepsilon\left(\displaystyle \frac{y}{|y|} \right)\displaystyle
\frac{y_i y_j}{|y|^2},\ \  \theta=\displaystyle \frac{y}{|y|} \in S^{n-1}.$$ Note that the function
$g_\varepsilon(\theta)$ is supported in $B_\epsilon(x_0)$ and $B_\epsilon(-x_0)$, where
$$B_\epsilon(x_0)=\{x\in S^{n-1}:|x-x_0|<\varepsilon \}.$$  The volume of these balls is of the order $\varepsilon^{n-1}$.

Now we want to show that $\displaystyle\frac{\partial^2}{\partial x_i \partial x_j}\left(
f_\varepsilon(\theta) r^{-n+1}\right) ^\wedge(x)$ can be made as small as desirable if
$\varepsilon$ is small. First let us show that $\displaystyle \frac{\partial^2}{\partial
y_i^2}\left( g_\varepsilon\left( \frac{y}{|y|}\right)\right)$ and therefore $\triangle
\left(g_\varepsilon(\theta) |y|^{-n+3} \right)(\xi)$ is of the order $\varepsilon^{-2}$.

To prove this, consider the function
 $$ f(x)=\left\{
\begin{array}{lll}
2 \displaystyle{e^{-\frac{|x-x_0|^2}{1-|x-x_0|^2}}}& \textrm{if } |x-x_0|<1 \\
2 \displaystyle{e^{-\frac{|x+x_0|^2}{1-|x+x_0|^2}}}& \textrm{if } |x+x_0|<1 \\
0 & \textrm{otherwise, }\end{array} \right.  $$ where $x\in S^{n-1}$.

The function $f$ is infinitely differentiable, so all its derivatives are bounded. Also
$f\left(\displaystyle\frac{x/|x|-x_0}{\varepsilon}+x_0\right)= f_{\varepsilon}(x/|x|)$. Therefore,

\begin{eqnarray*}
\left|\frac{\partial^k}{\partial x_i^k}f_\varepsilon(x/|x|) \right| =
\left|\frac{\partial^k}{\partial x_i^k} f\left(\displaystyle \frac{x/|x|-x_0} {\varepsilon}+x_0
\right ) \right| \le C_k \varepsilon^{-k},
\end{eqnarray*}
where $ C_k$ depends on $k$ but not on $x$.

The same is true for the derivatives of $g_\varepsilon(\theta)=f_\varepsilon\left(\displaystyle
\frac{y}{|y|} \right)\displaystyle \frac{y_i^2}{|y|^2}$, i.e.

$$\left|\frac{\partial^k}{\partial y_i^k}\ g_\varepsilon\left( \frac{y}{|y|}\right) \right|\le \tilde{C_k} \varepsilon^{-k}.$$

Therefore,
$$\left |\int_{S^{n-1}\cap (x/|x|_2)^\bot}\triangle \left(g_\varepsilon(\theta) |y|^{-n+3}
\right)(\xi) d\xi \right| \le $$ $$\le \tilde{C}_2\sup \left|\triangle \left(g_\varepsilon(\theta)
|y|^{-n+3} \right)\right| \cdot \int_{\left(B_\varepsilon(x_0)\cup B_\varepsilon(-x_0)\right)\cap
(x/|x|_2)^\bot}d\xi$$ $$= O(\varepsilon^{-2}\varepsilon^{n-2})=O(\varepsilon^{n-4}),$$
since the volume of the balls $B_\varepsilon(x_0)$ and $ B_\varepsilon(-x_0)$ is of the order
$\varepsilon^{n-1}$ and the volume of their intersection with a hyperplane is of the order
$\varepsilon^{n-2}$.

So $\displaystyle\frac{\partial^2}{\partial x_i  \partial x_j}\left(
f_\varepsilon(\theta)r^{-n+1}\right) ^\wedge(x)=O(\varepsilon^{n-4})$ and for $n>4$ it can be made
as small as desirable uniformly with respect to $x \in S^{n-1}$.

Using the same argument, we prove that $\left(f_\varepsilon(\theta)r^{-n+1}\right)
 ^\wedge(x)=O(\varepsilon^{n-2})$.

To get the estimate for $\displaystyle\frac{\partial}{\partial x_i}\left(
f_\varepsilon(\theta)r^{-n+1}\right) ^\wedge(x)$ we again use Lemma \ref{Lem:Kold} and the
connection between the Fourier transform and differentiation. Take $\alpha$ close to $1$. Then, if
$g_\varepsilon(\theta)=f_\varepsilon\left(\displaystyle \frac{y}{|y|} \right)\displaystyle
\frac{y_i}{|y|},\ \  \theta=\displaystyle \frac{y}{|y|} \in S^{n-1},$

$$\displaystyle\frac{\partial}{\partial x_i}\left( f_\varepsilon(\theta)r^{-n+\alpha+1}\right)
^\wedge(x)= -\left( g_\varepsilon(\theta) |y|^{-n+\alpha+2}  \right)^\wedge (x) $$
\begin{equation}\label{eqn:limit}
=\frac{-\pi |x|_2^{-2} \int_{S^{n-1}} |(x,\xi)|^{1-\alpha}\Delta
(g_\varepsilon(\theta)|y|^{-n+\alpha+1})(\xi)d\xi}{2 \Gamma(2-\alpha)\sin\frac{\pi(1-\alpha)}{2}}.
\end{equation}
When $\alpha$ approaches $1$, the numerator and denominator in the right hand side approach zero.
Indeed, let us show that the limit of the numerator is zero:
$$\lim_{\alpha\to 1}\int_{S^{n-1}} |(x,\xi)|^{1-\alpha}\Delta
(g_\varepsilon(\theta)|y|^{-n+\alpha+1})(\xi)d\xi=\int_{S^{n-1}}\Delta
(g_\varepsilon(\theta) |y|^{-n+2})(\xi)d\xi.$$ Recall the relation between the
spherical Laplacian $\Delta_S$ and Euclidean Laplacian $\Delta$ (see, for example,
\cite[p.7]{Gr}): if $f$ is a homogeneous function of degree $m$, then on the sphere
$$\Delta_S f=\Delta f -m(m+n-2)f.$$
Since $g_\varepsilon(\theta)|y|^{-n+2}$ has degree of homogeneity $-n+2$, the previous formula
implies $\Delta (g_\varepsilon(\theta)|y|^{-n+2})(\xi)= \Delta_S(g_\varepsilon
(\theta)|y|^{-n+2})(\xi)$. Due to the fact that $\Delta_S$ is a self-adjoint operator,
\cite[p.7]{Gr}, we have
$$\int_{S^{n-1}}\Delta_S (g_\varepsilon(\theta)|y|^{-n+2})(\xi)d\xi=0.$$

Now to compute the limit of (\ref{eqn:limit}) as $\alpha \to 0$, apply l'Hopital's rule:
\begin{eqnarray*}
\displaystyle\frac{\partial}{\partial x_i}\left( f_\varepsilon(\theta)r^{-n+2}\right) ^\wedge(x)&=&
|x|_2^{-2} \int_{S^{n-1}} \ln|(x,\xi)|\Delta (g_\varepsilon(\theta)|y|^{-n+2})(\xi)d\xi.
\end{eqnarray*}
Recall that the function $g_\varepsilon(\theta) $ is supported in the balls $B_\varepsilon(x_0)$
and $B_\varepsilon(-x_0)$. Then
$$\left|\int_{S^{n-1}}\ln|(x,\xi)|\Delta (g_\varepsilon(\theta)|y|^{-n+2})(\xi)d\xi\right|\le \sup
\left|\Delta (g_\varepsilon(\theta)|y|^{-n+2})(\xi)\right|\times$$  $$\times
\int_{B_\varepsilon(x_0)\cup B_\varepsilon(-x_0)}\left|\ln|(x,\xi)|\right|d\xi.$$

Now we want to estimate the latter integral. Note that it is enough to estimate just
$\int_{B_\varepsilon(x_0)}\left|\ln|(x,\xi)|\right|d\xi$. Consider two cases. First, suppose that
$x $ is not perpendicular to any $y \in B_{2\varepsilon}(x_0)$. In this case one can check that
$|(x,\xi)|>\varepsilon/2 $ and therefore
$$\int_{B_\varepsilon(x_0)}\left|\ln|(x,\xi)|\right|d\xi\le \left|\ln (\varepsilon/2) \right| \cdot
\mbox{vol}(B_\varepsilon(x_0))=O(|\ln\varepsilon|\cdot\varepsilon^{n-1})$$ since the volume of the
ball $B_\varepsilon(x_0)$  is of the order $\varepsilon^{n-1}$.

In the second case there exists  $y\in B_{2\varepsilon}(x_0)$ such that $x \bot y$. Consider the
ball $B_{4\varepsilon}(y)$. Clearly, $B_{2\varepsilon}(x_0)\subset B_{4\varepsilon}(y)$, therefore
$$\int_{B_\varepsilon(x_0)}\left|\ln|(x,\xi)|\right|d\xi\le \int_{B_{4\varepsilon}(y)}\left|\ln|(x,\xi)|\right|d\xi.$$

Let us make a change of coordinates from $\xi \in  S^{n-1}$ to $\zeta \in S^{n-2} $ and $t \in
[-\pi/2,\pi/2]$ such that $\xi=y\sqrt{1-t^2}+\zeta t$. The Jacobian is equal to $\displaystyle
\frac{t^{n-2}}{\sqrt{1-t^2}}$. If $\xi \in B_{4\varepsilon}(y)$ then $t \in
[0,4\varepsilon\sqrt{1-4\varepsilon^2}]\subset[0,4\varepsilon]$. Using the fact that $x$ is
perpendicular to $y$, we get

$$ \int_{B_{4\varepsilon}(y)}\left|\ln|(x,\xi)|\right|d\xi\le \int_{S^{n-2}}\int_{0}^{4\epsilon}\left|
 \ln|(x,\zeta)t|\right|\frac{t^{n-2}}{\sqrt{1-t^2}} dt d\zeta$$
$$=\int_{S^{n-2}}\left(\int_{0}^{4\epsilon}\left| \ln|(x,\zeta)|\right|\frac{t^{n-2}}
{\sqrt{1-t^2}}dt+\int_{0}^{4\epsilon}\left| \ln t\right|\frac{t^{n-2}}{\sqrt{1-t^2}}dt\right) d\xi
$$
$$\le 2 \int_{S^{n-2}}\left(\int_{0}^{4\epsilon}\left| \ln|(x,\zeta)|\right| t^{n-2} dt+
\int_{0}^{4\epsilon}\left| \ln t\right|t^{n-2}dt\right) d\xi $$ for a small $\varepsilon$. The
first integral can be estimated in the following way:

$$\int_{S^{n-2}}\left| \ln|(x,\zeta)|\right|d\xi \int_{0}^{4\epsilon}t^{n-2} dt = O(\varepsilon^{n-1}).$$

Using integration by parts in the second integral, we get that

$$\int_{0}^{4\epsilon}\left| \ln t\right|t^{n-2}dt =\frac{1}{n-1} \ln t \cdot t^{n-1}\Big|_0^{4\varepsilon}-
\frac{1}{n-1}\int_{0}^{4\epsilon}t^{n-2}dt$$
$$ =O(|\ln\varepsilon|\cdot\varepsilon^{n-1}).$$

Therefore, $$\int_{B_\varepsilon(x_0)\cup B_\varepsilon(-x_0)}\left|\ln|(x,\xi)|\right|d\xi =
O(|\ln\varepsilon|\cdot\varepsilon^{n-1})$$ and $$\displaystyle\frac{\partial}{\partial x_i}\left(
f_\varepsilon(\theta)r^{-n+\alpha+1}\right)^\wedge(x)= O(|\ln\varepsilon|\cdot\varepsilon^{n-3}).$$



 \qed

\begin {Lem}\label{K-convex}
If $n\ge5$, the body $K$ is convex for small enough $\varepsilon$.
\end{Lem}

\noindent{\bf Proof.} Let $|\cdot|_2$ be the Euclidean norm. By \cite[p.363]{GS} the Fourier
transform of $|x|_2^q$, $q \in (-n,0)$ equals
$$(|x|_2^q)^\wedge(t)=2^{q+n}\pi^{n/2}\frac{\Gamma(\frac{q+n}{2})} {\Gamma(\frac{-q}{2})}|t|_2^{-n-q} .$$

Using this formula and the definition of the body $K$
\begin{eqnarray}\label{bodyK}
\|x\|_K^{-1}&=&\left( (1- f_\varepsilon(\theta))r^{-n+1}\right)^\wedge(x) \nonumber \\
&=& C_n|x|_2^{-1} - \left( f_\varepsilon(\theta)r^{-n+1}\right)^\wedge(x),
\end{eqnarray}
where $C_n=\displaystyle\frac{2\pi^{(n+1)/2}}{\Gamma(\frac{n-1}{2})}$.

Let $K_W$ be the section of $K$ by a $2$-dimensional central plane $W$ with an orthonormal basis
$\xi_1$, $\xi_2$. So, if $x \in W \cap S^{n-1}$, then $x = \xi_{1}\cos\phi+\xi_{2}\sin\phi$, $\phi
\in [0,2\pi]$. To show that $K$ is convex, it is enough to show that $K_W$ is convex for any $W$.

Consider a function
$$\rho_{W}(\phi)= \left( f_\varepsilon(\theta)r^{-n+1}\right)^\wedge(\xi_{1}\cos\phi+\xi_{2}\sin\phi).$$
By the definition of $K$, the radial function of $K_W$ is given by
$$\rho(\phi)=C_n - \rho_W(\phi). $$
To prove that $K_W$ is convex, we need to show that for small $\varepsilon$

$$J(W, \varepsilon, \phi)=2(\rho')^2-\rho''\rho+\rho^2 > 0 $$
for every $W$ and $\phi$, see \cite[p.25]{Ga2}.

Computing the derivatives,
\begin{eqnarray*}
\rho'(\phi)&=& -\frac{d}{d \phi}\left(\rho_W(\phi) \right) \\
\rho''(\phi)&=& -\frac{d^2}{d ^2\phi}\left( \rho_W(\phi) \right).
\end{eqnarray*}

To estimate $\rho'$ and $\rho''$, consider
\begin{eqnarray*}
\left|\frac{d}{d \phi}\left( \rho_W(\phi) \right) \right| &=&
\left|\sum_{i=1}^n\frac{\partial}{\partial x_i}\left(
f_\varepsilon(\theta)r^{-n+1}\right)^\wedge(x)\right|\left| \frac{d x_i}{d \phi}
\right|\\
 &\le & 2\left|\sum_{i=1}^n\frac{\partial}{\partial
 x_i}\left(f_\varepsilon(\theta)r^{-n+1}\right)^\wedge(x)\right|,
\end{eqnarray*}
since $x = \xi_{1}\cos\phi+\xi_{2}\sin\phi$ and
$$ \left| \frac{d x_i}{d \phi} \right|= \left|-\xi_{1,i}\sin\phi+\xi_{2,i}\cos\phi
\right|\le 2. $$

Similarly, $$\left|\frac{d^2}{d \phi^2}\left( \rho_W(\phi) \right) \right| \le
4\left|\sum_{i,j=1}^n \frac{\partial^2} {\partial x_i \partial
x_j}\left(f_\varepsilon(\theta)r^{-n+1}\right)^\wedge(x)\right|.$$

By Lemma \ref{Lem:est} we have $\rho'=O(|\ln \varepsilon|\varepsilon^{n-3})$ and
$\rho''=O(\varepsilon^{n-4})$. Since these estimates are uniform with respect to $\phi$ and $W$, it
follows that, for small enough $\varepsilon$, $J(W, \varepsilon, \phi)>0$ for every $\phi$ and $W$.

\qed

In \cite[Lemma 1]{K-1}  the following  was proved. Let $f$ be an even continuous homogeneous
function of degree $-n+1$ on $\mathbb{R}^n \setminus \{0\}$. Then $\hat{f}$ is a continuous
function outside of the origin and for every $\xi \in S^{n-1}$,
\begin{equation}\label{Fourier-sphere}
\hat{f}(\xi)=\pi\int_{S^{n-1}\cap\{(\theta,\xi)=0\}}f(\theta)d\theta.
\end{equation}

Using this formula,
\begin{equation}\label{int-formula}
\left( (1- f_\varepsilon(\theta))r^{-n+1}\right)^\wedge(x)=\pi\int_{S^{n-1}\cap
\{(x,\theta)=0\}}(1- f_\varepsilon(\theta))d\theta,
\end{equation}

\begin{Lem}\label{V-int}
For every $(n-1)$-dimensional subspace $V$ of $\mathbb{R}^n$, the body $K \cap V$ is an
intersection body.
\end{Lem}

\noindent{\bf Proof.} Fix an $(n-1)$-dimensional subspace $V$ of $\mathbb{R}^n$. By the definition
(\ref{Def:K}) of the body $K$ and formula (\ref{int-formula})

\begin{eqnarray*}
\frac{1}{\pi}\|x\|^{-1}_{K}&=&\int_{S^{n-1}\cap \{(x,\theta)=0\}} (1-f_\varepsilon(\theta))d\theta.
\end{eqnarray*}
In particular, for $x \in V$,
\begin{eqnarray*}
\frac{1}{\pi}\|x\|^{-1}_{K\cap V}&=&\int_{S^{n-1}\cap \{(x,\theta)=0, x \in V  \}}
(1-f_\varepsilon(\theta))d\theta,
\end{eqnarray*}
where $\{(x,\theta)=0, x \in V\}$ is the hyperplane of all $\theta$ that are perpendicular to $x$
for a fixed $x\in V$.


Let us change coordinates from $\theta \in S^{n-2}=S^{n-1}\cap \{(x,\theta)=0, x \in V  \}$ to
$\phi \in [0,\pi]$ and $\xi \in S^{n-3}=S^{n-2}\cap \{(x,\theta)=0, x \in V  \}$. The Jacobian is
equal to $(\sin \phi)^{n-3}$. So

$$\frac{1}{\pi}\|x\|^{-1}_{K\cap V}=\int_{S^{n-2}\cap \{(x,\theta)=0, x \in V\}}\left(\int_0^{\pi}
(1-f_\varepsilon(\xi,\phi))(\sin \phi)^{n-3} d\phi \right)d\xi $$ $$= \int_{S^{n-2}\cap
\{(x,\theta)=0, x \in V\}}\left(\int_0^{\pi}(\sin \phi)^{n-3} d\phi-
\int_0^{\pi}f_\varepsilon(\xi,\phi)(\sin \phi)^{n-3} d\phi\right)d\xi.$$

Taking the Fourier transform of both sides, as functions of the variable $x \in V$, and using
(\ref{Fourier-sphere})
\begin{eqnarray*}
\left(\|x\|^{-1}_{K\cap V}\right)^\wedge(\theta)&=&\pi\int_0^{\pi}(\sin \phi)^{n-3} d\phi-
\pi\int_0^{\pi}f_\varepsilon(\theta)(\sin \phi)^{n-3} d\phi
\end{eqnarray*}
for $\theta \in S^{n-2}$.

By the definition, $f_\varepsilon$ is non-zero only in an $\varepsilon$-neighborhood of $x_0$, so
there exists a set $R_\varepsilon \subset [0,\pi]$ such that $f_\varepsilon(\xi,\phi)=0$ for $\phi
\in [0,\pi]\setminus R_\varepsilon$. Also, $|f_\varepsilon|\le 2$ and $|R_\varepsilon|$, the length
of the one-dimensional set $R_\varepsilon$, is of the order $\varepsilon$. Therefore,
\begin{eqnarray*}
\int_0^{\pi}f_\varepsilon(\xi,\phi)(\sin \phi)^{n-3} d\phi&=&\int_{R_\varepsilon} f_\varepsilon
(\xi,\phi)(\sin \phi)^{n-3} d\phi\\
&\le& 2|R_\varepsilon|=C \varepsilon,
\end{eqnarray*}
where $C$ does not depend of the choice of $V$.

Since $\int_0^{\pi}(\sin \phi)^{n-3} d\phi$ is equal to some positive number,
$$\int_0^{\pi}(\sin
\phi)^{n-3} d\phi- \int_0^{\pi}f_\varepsilon(\theta)(\sin \phi)^{n-3} d\phi>0$$ for
a sufficiently small $\varepsilon$, which means we can find $\varepsilon$ small
enough that $\left(\|x\|^{-1}_{K\cap V}\right)^\wedge(\theta)>0$ for all $\theta$.
Therefore by \cite[Thm.1]{K1}, for small enough $\varepsilon$, for every
$(n-1)$-dimensional subspace $V$ of $\mathbb{R}^n$, the body $K \cap V$ is an
intersection body. \qed

Now Theorem \ref{Thm:main} follows from Lemmas \ref{K-int}, \ref{K-convex} and \ref{V-int}.

\noindent{\bf Remark 1.} By a result of Neyman \cite{N}, for any finite system of equations and
inequalities involving the norms, there exists a subspace of $L_p$, $0<p<2$ that does not satisfy
this system. Since, by \cite{K1}, the unit ball of every finite dimensional subspace of $L_p$,
$0<p<2$ is an intersection body, we conclude that the class of intersection bodies cannot be
characterized by a finite number of equations or inequalities.

\noindent{\bf Remark 2.} For $n=5$ the fact that there exists a non-intersection
body whose central sections are intersection bodies follows from the solution of the
Busemann-Petty problem since every four-dimensional symmetric convex body is an
intersection body, see \cite{GKS}, \cite{Zh2}.

{\bf Acknowledgments.} The author is thankful to A.Koldobsky and V.Yaskin for their
valuable help in preparation of the paper.


\end{document}